\documentclass[envcountsect]{svmult}
\pdfoutput=1
\usepackage{graphicx}
\usepackage{amssymb, amsmath}
\usepackage[bookmarks=false]{hyperref}
\usepackage{caption}                                
\allowdisplaybreaks[1] 

\spnewtheorem{Theorem}[theorem]{Theorem}{\bfseries}{\itshape}
\spnewtheorem{Lemma}[theorem]{Lemma}{\bfseries}{\itshape}
\spnewtheorem{Corollary}[theorem]{Corollary}{\bfseries}{\itshape}
\spnewtheorem{Proposition}[theorem]{Proposition}{\bfseries}{\itshape}
\spnewtheorem{Remark}[theorem]{Remark}{\itshape}{\rmfamily}
\spnewtheorem{Definition}[theorem]{Definition}{\bfseries}{\itshape}
\spnewtheorem{Example}{Example}{\bfseries}{\itshape}
\spnewtheorem{Question}{Question}{\bfseries}{\itshape}
\spnewtheorem{Answer}[theorem]{Answer}{\bfseries}{\itshape}

\newcommand{\Q}{\mathbb{Q}} 
\newcommand{\Z}{\mathbb{Z}} 
\def\S{{\cal S}}

\smartqed

\begin{document} 

\title*{Partition-theoretic formulas for arithmetic densities}
\titlerunning{Partition-theoretic formulas for arithmetic densities}

\author{Ken Ono, Robert Schneider, and Ian Wagner \thanks{This research was partially supported by NSF grant
     DMS 1601306. The first author thanks the support of the Asa Griggs Candler Fund, and the second author thanks the support of an Emory University Woodruff Fellowship.}}
\institute{Dept. Mathematics and Computer Science, Emory University, Atlanta, Georgia 30322\\
\email{ono@mathcs.emory.edu}, \
\email{robert.schneider@emory.edu}, \newline
\email{iwagner@emory.edu}
}
\motto{In celebration of Krishnaswami Alladi's 60th birthday.}
\maketitle


\abstract{ 
If $\gcd(r,t)=1$, then Alladi proved the M\"obius sum identity
$$-\sum_{\substack{ n \geq 2 \\ p_{\rm{min}}(n) \equiv r \pmod{t}}} \mu(n)n^{-1}= \frac{1}{\varphi(t)}.
$$ 
Here $p_{\rm{min}}(n)$ is the smallest prime divisor of $n$.  The right-hand side represents the proportion of primes in a fixed arithmetic progression modulo $t$.   Locus generalized this to Chebotarev densities for Galois extensions.  Answering a question of Alladi, we obtain analogs of these results to arithmetic densities of subsets of positive integers using $q$-series and integer partitions.  For suitable subsets $\S$ of the positive integers with density $d_{\S}$, we prove that 
\[- \lim_{q \to 1} \sum_{\substack{ \lambda \in \mathcal{P} \\ \rm{sm}(\lambda) \in \S}} \mu_{\mathcal{P}} (\lambda)q^{\vert \lambda \vert} =  d_{\S},\]
where the sum is taken over integer partitions $\lambda$, $\mu_{\mathcal{P}}(\lambda)$ is a partition-theoretic M\"obius function,   
 $\vert \lambda \vert$ is the size of partition $\lambda$, and $\rm{sm}(\lambda)$ is the smallest part of $\lambda$.   In particular, we obtain partition-theoretic formulas for even powers  of $\pi$ when considering power-free integers.
\keywords{Arithmetic densities, partition-theoretic M\"obius function}\\
{\bf 2010 Mathematics Subject Classification.} {05A17, 11P82}
\\[12pt]
}

\renewcommand{\theequation}{\arabic{section}.\arabic{equation}}

\section{Introduction and Statement of Results}

The  fact 
$\sum_{n=1}^{\infty} \mu(n)/n= 0,$
where $\mu(n)$ is the  M\"{o}bius function, can be reformulated as
\[-\sum_{n=2}^{\infty} \frac{\mu(n)}{n} = 1.\]
For notational convenience  define $\mu^{*}(n) := - \mu(n)$.  This equation above can be interpreted as the statement that one-hundred percent 
of integers $n \geq 2
$ are divisible by at least one prime. This idea was used by Alladi  \cite{A} to prove that if $\gcd(r, t) = 1$, then
\begin{equation}\label{AlladiSum}
\sum_{\substack{ n \geq 2 \\ p_{\rm{min}}(n) \equiv r \pmod{t}}} \frac{\mu^{*}(n)}{n} = \frac{1}{\varphi(t)}. 
\end{equation}
Here $\varphi(t)$ is Euler's phi function, and $p_{\rm{min}}(n)$ is the smallest prime factor of $n$.  

Alladi has asked \cite{Alladi_lecture} for a partition-theoretic generalization of this result.  We answer his question by obtaining an analog of a generalization that was recently obtained by Locus \cite{L}.  Locus began by interpreting Alladi's theorem as a device for computing densities of primes in arithmetic progressions.  She generalized this idea, and proved analogous formulas for the Chebotarev densities of Frobenius elements in unions of conjugacy classes of Galois extensions. 

We recall Locus's result.  Let $S$ be a subset of primes with Dirichlet density, and define
\begin{equation}\label{maddie}
\mathfrak{F}_{S}(s) := \sum_{\substack{n \geq 2 \\ p_{\rm{min}}(n) \in S}} \frac{\mu^{*}(n)}{n^{s}}.  
\end{equation}
Suppose $K$ is a finite Galois extension of $\Q$ and $p$ is an unramified prime in $K$.  Define
\[\left[ \frac{K/\Q}{p} \right] := \left\{ \left[ \frac{K/\Q}{\mathfrak{p}} \right] : \mathfrak{p} \subseteq \mathcal{O}_{K} \ \rm{is \  a \  prime \  ideal \  above} \ \it{p} \right\},\]
where $\left[ \frac{K/\Q}{\mathfrak{p}} \right]$ is the Artin symbol (for example, see Chapter 8 of \cite{Artin_ref_needed}), and $\mathcal{O}_{K}$ is the ring of integers of $K$.  It is well known that $\left[ \frac{K/\Q}{p} \right]$ is a conjugacy class $C$ in $G= \rm{Gal}( \emph{K}/\Q)$.  If we let
\begin{equation}\label{OnePointThree}
S_{C} := \left\{ p \ \rm{prime}: \left[ \frac{\emph{K}/\Q}{\emph{p}} \right] = C \right\},
\end{equation}
then Locus proved (see Theorem $1$ of \cite{L}) that 
$$
\mathfrak{F}_{S_{C}}(1) = \frac{ \# C}{\# G}. 
$$

\begin{remark}
Alladi's formula (\ref{AlladiSum}) is the cyclotomic case of Locus' Theorem.
\end{remark}

We now turn to Alladi's question concerning a partition-theoretic analog.   A \textit{partition} is a finite non-increasing sequence of positive integers, say $$\lambda = (\lambda_{1}, \lambda_{2},..., \lambda_{\ell(\lambda)}),$$
 where $\ell(\lambda)$ denotes the length of $\lambda$, i.e. the number of parts.  The size of $\lambda$ is $\vert \lambda \vert := \lambda_{1} + \lambda_{2} + \cdot \cdot \cdot + \lambda_{\ell(\lambda)}$, i.e., the number being partitioned.  Furthermore, we let $\rm{sm}(\lambda) := \lambda_{\ell(\lambda)}$ denote the smallest part of $\lambda$ (resp. $\rm{lg}(\lambda) := \lambda_{1}$ the largest part of $\lambda$).  
We require the partition-theoretic M\"{o}bius function
\begin{equation}\label{OnePointFour}
\mu_{\mathcal{P}}(\lambda) := \begin{cases} 0 & \rm{if} \ \lambda \ \rm{has \ repeated \ parts} \\
(-1)^{\ell(\lambda)} & \rm{otherwise} \end{cases}  
\end{equation}
examined in \cite{Schneider_arithmetic}, where $\mathcal{P}$ denotes the set of integer partitions.  Notice that $\mu_{\mathcal{P}}(\lambda) = 0$ if $\lambda$ has any repeated parts, which is analogous to the vanishing of $\mu(n)$ for integers $n$ which are not square-free.  In particular, the parts in partition $\lambda$ play the role of prime divisors of $n$ in this analogy.  We define $\mu_{\mathcal{P}}^{*}(\lambda):= - \mu_{\mathcal{P}}(\lambda)$ as in Locus' theorem, for aesthetic reasons.  

The table below offers a complete description of the objects which are related with respect to this analogy.  However, it is worthwhile to first compare the generating functions for $\mu(n)$ and $\mu_{\mathcal{P}}(\lambda)$.  Using the Euler product for the Riemann zeta function, it is well known that the Dirichlet generating function for $\mu(n)$ is
\begin{equation}\label{OnePointFive}
\frac{1}{\zeta(s)} = \prod_{p \ \rm{prime}} \left( 1 - \frac{1}{p^s} \right) = \sum_{m=1}^{\infty} \mu(m)m^{-s}. 
\end{equation}
To obtain the generating function for $\mu_{\mathcal{P}}(\lambda)$, we recall the $q$-Pochhammer symbol 
\begin{equation}\label{OnePointSix}
(a;q)_{n} := \prod_{m=0}^{n-1} (1-aq^{m}).
\end{equation}
For $|q|<1$, the $q$-Pochhammer symbol with $n = \infty$ is naturally defined by 
\begin{equation}\label{OnePointSeven}
(a;q)_{\infty} := \lim_{n\to \infty}(a;q)_{n}. 
\end{equation}
Using this symbol, it is clear that the generating function for $\mu_{\mathcal{P}}(\lambda)$ is  
$$
(q;q)_{\infty} = \prod_{n=1}^{\infty} (1 - q^n) = \sum_{\lambda} \mu_{\mathcal{P}}(\lambda) q^{\vert \lambda \vert}. 
$$
By comparing the generating functions for $\mu(n)$ and $\mu_{\cal{P}}(\lambda)$,  we see that prime factors and integer parts
of partitions are natural analogs of each other. The following table offers the identifications that make up this analogy. 

\vspace{0.5cm}

\begin{center}
\begin{tabular}{ | c | c | }
\hline
Natural number $m$ & Partition $\lambda$\\ \hline
Prime factors of $m$ & Parts of $\lambda$ \\ \hline
Square-free integers & Partitions into distinct parts \\ \hline
$\mu(m)$ & $\mu_{\mathcal{P}}(\lambda)$ \\ \hline
$p_{\rm{min}}(m)$ & $\rm{sm}(\lambda)$ \\ \hline
$p_{\rm{max}}(m)$ & $\rm{lg}(\lambda)$ \\ \hline
$m^{-s}$ & $q^{\vert \lambda \vert}$ \\ \hline
${\zeta(s)^{-1}}$ & ${(q;q)_{\infty}}$ \\  \hline
$s=1$ & $q\rightarrow 1$\\
\hline
\end{tabular}
\end{center}
\vspace{0.5cm}

\begin{remark}There are further  analogies between multiplicative number theory and the theory of integer partitions. 
For instance, in  \cite{AlladiErdos} Alladi and Erd\H{o}s exploited a bijection between prime factorizations of integers and partitions into prime parts, to study an interesting arithmetic function; recently, the first two authors have shown that many theorems in multiplicative number theory are special cases of much more general partition-theoretic phenomena \cite{OnoRolenS,Schneider_arithmetic}.
\end{remark}

Suppose that
$\S$ is a subset of the positive integers with arithmetic density
\[\lim_{X \to \infty} \frac{\# \{ n \in \S : n \leq X \}}{X} = d_{\S}.\]
The partition-theoretic counterpart to (\ref{maddie}) is 
\begin{equation}\label{OnePointNine}
F_{\mathcal{S}}(q) := \sum_{\substack{ \lambda \in \mathcal{P} \\ \rm{sm}(\lambda) \in \mathcal{S}}} \mu_{\mathcal{P}}^{*}(\lambda) q^{\vert \lambda \vert}. 
\end{equation}
To state our results, we define 
\begin{equation}\label{OnePointTen}
S_{r,t} := \{n \in \Z^{+} : n \equiv r \pmod{t} \}. 
\end{equation}
These sets are simply the positive integers in an arithmetic progression $r$ modulo $t$. 

Our first result concerns the case where $t=2$. Obviously, the arithmetic densities of $\S_{1,2}$ and $\S_{2,2}$ are both $\frac{1}{2}$. The theorem below offers a formula illustrating these densities and also offers curious lacunary $q$-series
identities expressed in terms of theta series.

\begin{Theorem} \label{12}
Assume the notation above.
\begin{enumerate}
\item The following $q$-series identities are true: 
\[F_{\S_{1,2}}(q) = \sum_{n=1}^{\infty} (-1)^{n+1} q^{n^2},\]
\[F_{\S_{2,2}}(q) = 1+ \sum_{n=1}^{\infty} (-1)^{n} q^{n^2} - \sum_{m= -\infty}^{\infty} (-1)^{m}q^{\frac{m(3m-1)}{2}}.\]
\item We have that
\[\lim_{q \to 1} F_{\S_{1,2}}(q) = \lim_{q \to 1} F_{\S_{2,2}}(q) = \frac{1}{2}.\]
\end{enumerate}
\end{Theorem}

\begin{remark}
An identity which implies the first part of Theorem~\ref{12} was obtained earlier by Andrews and Stenger \cite{AndrewsStenger}. This identity and a generalization is also discussed in Examples 2 and 3 on p. 156-157 of \cite{And}.
\end{remark}

\begin{remark}
The limits in Theorem \ref{12} are understood as $q$ tends to $1$ from within the unit disk.
\end{remark}

\begin{Example}
For complex $z$ in the upper-half of the complex plane, let $q(z) := \rm{exp} \left( - \frac{2 \pi i}{\emph{z}} \right)$.  Therefore, if $z \to 1$ in the upper-half plane, then $q(z) \to 1$ in the unit disk.  The table below displays a set of such $z$ beginning to approach $1$ and the corresponding values of $F_{\S_{1,2}}(q(z))$.

\begin{center}
\begin{tabular}{ | c | c | }
\hline
$z$ & $F_{\S_{1,2}}(q(z))$\\ \hline
$1+.10i$ & $0.458233...$ \\ \hline
$1+.09i$ & $0.471737...$ \\ \hline
$1 + .08i$ & $0.482784...$ \\ \hline
$1+.07i$ & $0.491003...$ \\ \hline
$1+.06i$ & $0.496296...$ \\ \hline
$1+.05i$ & $0.498998...$ \\ \hline
$1+.04i$ & $0.499919...$ \\ \hline
$1+.03i$ & $0.500048...$ \\ \hline
$1+.02i$ & $0.500024...$ \\ \hline
$1+.01i$ & $0.500006...$ \\ 
\hline
\end{tabular}
\end{center}
\vspace{0.5cm}
\end{Example}

The first claim in Theorem \ref{12}  offers an immediate combinatorial interpretation.  Let $D^{+}_{\rm{even}}(n)$ denote the number of partitions of $n$ into an even number of distinct parts with smallest part even, and let $D^{+}_{\rm{odd}}(n)$ denote the number of partitions of $n$ into an even number of distinct parts with smallest part odd.  Similarly, let $D^{-}_{\rm{even}}(n)$ denote the number of partitions of $n$ into an odd number of distinct parts with smallest part even, and let $D^{-}_{\rm{odd}}(n)$ denote the number of partitions of $n$ into an odd number of distinct parts with smallest part odd. 
To make this precise, for integers $k$ let $\omega(k):=\frac{k(3k-1)}{2}$ 
be the index $k$ {\it pentagonal number}.

\begin{Corollary} \label{cor}
Assume the notation above.
\begin{enumerate}
\item For partitions into distinct parts whose smallest part is odd, we have
\begin{displaymath}
D_{\rm{odd}}^{+}(n) - D_{\rm{odd}}^{-}(n) \\
= \begin{cases} 0 & {\rm{if}} \  n \ \rm{is \ not \ a \ square} \\
1 & {\rm{if}} \ n \ \rm{is \ an \ even \ square} \\
-1 & {\rm{if}} \ n \ \rm{is \ an \ odd \ square}. \end{cases}
\end{displaymath}

\item For partitions into distinct parts whose smallest part is even, we have
\begin{displaymath}
\begin{split}
&D_{\rm{even}}^{+}(n) - D_{\rm{even}}^{-}(n)\\
&\ \ \ \ \  = \begin{cases}
-1 & {\rm{if}} \ n \ \rm{is \ an \ even \ square \ and \ not \ a \ pentagonal \ number} \\
1 & {\rm{if}} \ n \ \rm{is \ an \ odd \ square \ and \ not \ a \ pentagonal \ number} \\
1 & {\rm{if}} \ n \ \rm{is \ an \ even\ index \ pentagonal \ number \ and \ not \ a \ square} \\
-1 & {\rm{if}} \ n \ \rm{is \ an \ odd \ index \ pentagonal \ number \ and \ not \ a \ square} \\
0 & {\rm{otherwise.}} \\
 \end{cases}
 \end{split}
 \end{displaymath}
\end{enumerate}
\end{Corollary}

\begin{Question}
It would be interesting to obtain a combinatorial proof of Corollary \ref{cor} which is in the spirit of Franklin's proof of Euler's Pentagonal Number Theorem (see pages 10-11 of \cite{And}).
\end{Question}

Our proof of Theorem \ref{12} makes use of the $q$-Binomial Theorem and some well-known $q$-series identities.  It is natural to ask whether such a relation holds for general sets $\S_{r,t}$.  The following theorem shows that Theorem \ref{12} is indeed a special case of a more general phenomenon.  

\begin{Theorem} \label{rt}
If  $0\leq r<t$ are integers and $\gcd(m,t)=1$, then we have that
$$\lim_{q \to \zeta} F_{\S_{r,t}}(q) = \frac{1}{t},
$$
where $\zeta$ is a primitive $m$th root of unity.
\end{Theorem}

\begin{remark}
The limits in Theorem \ref{rt} are understood as $q$ tends to $\zeta$ from within the unit disk.
\end{remark}

Obviously, these results hold for any set $\S$ of positive integers that is a finite union of arithmetic progressions.  It turns out that this theorem can also be used to compute arithmetic densities of more complicated sets arising systematically from a careful study of arithmetic progressions.  We focus on the sets of positive integers $\S^{(k)}_{\rm{fr}}$ which are $k$th power-free.  In particular, we have that \[\S^{(2)}_{\rm{fr}} = \{ 1, 2, 3, 5, 6, 7, 10, 11, 13, \dots \}.\]
It is well known that the arithmetic densities of these sets are given by
$$
\lim_{X\rightarrow +\infty}\frac{\# \left \{1\leq n\leq X \ : \ n\in \S_{\rm{fr}}^{(k)}\right \}}{X}=
\prod_{p \ {\text {\rm prime}}} \left (1-\frac{1}{p^k}\right)=\frac{1}{\zeta(k)}.
$$
To obtain partition-theoretic formulas for these densities, we first compute a partition-theoretic formula for the density of 
\begin{equation}
\S^{(k)}_{\rm{fr}}(N):= \{n \geq 1 \ : \ p^{k} \nmid n \ \rm{for \ every} \ \emph{p} \leq \emph{N} \}. 
\end{equation}

\begin{Theorem} \label{free}
If $k, N \geq 2$ are integers, then we have that
\[ \lim_{q \to 1} F_{\S^{(k)}_{\rm{fr}}(N)}(q) = \prod_{p \leq N \ \rm{prime}} \left(1 - \frac{1}{p^k} \right).\]
\end{Theorem}

The constants in Theorem \ref{free} are the arithmetic densities of positive integers that are not divisible by the $k$th power of any prime $p \leq N$, namely $\S^{(k)}_{\rm{fr}}(N)$.
Theorem \ref{free} can be used to calculate the arithmetic density of $\S^{(k)}_{\rm{fr}}$ by letting $N\rightarrow +\infty$.

\begin{Corollary} \label{free2}
If $k \geq 2$, then 
\[\lim_{q \to 1} F_{\S^{(k)}_{\rm{fr}}}(q) = \frac{1}{\zeta(k)}.\]
Furthermore, if $k \geq 2$ is even, then 
\[\lim_{q \to 1} F_{\S^{(k)}_{\rm{fr}}}(q) =(-1)^{\frac{k}{2}+1} \frac{k!}{B_{k}\cdot 2^{k-1}}\cdot \frac{1}{\pi^k},\]
where $B_{k}$ is the $k$th Bernoulli number.
\end{Corollary}

This paper is organized as follows.  In Section $2.1$ we discuss the $q$-Binomial Theorem, which will be an essential tool for our proofs, as well as a duality principle for partitions related to ideas of Alladi.  In Section $2.2$ we will use the $q$-Binomial Theorem to prove results related to Theorem \ref{rt}.  Section $3$ will contain the proofs of all of the theorems, and Section $4$ will contain some nice examples. 

\section{The $q$-Binomial Theorem and its consequences}

In this section we recall elementary $q$-series identities, and we offer convenient reformulations for
the functions $F_{\mathcal{S}}(q)$.

\subsection{Easy nuts and bolts}
Let us recall the classical $q$-Binomial Theorem (see \cite{And} for proof). 

\begin{Lemma}\label{lemma1}
For $a,z\in \mathbb C, |q|<1$ we have the identity
\begin{equation*}
\frac{(az;q)_{\infty}}{(z;q)_{\infty}}=\sum_{n=0}^{\infty}\frac{(a;q)_n}{(q;q)_n}\cdot z^n.
\end{equation*}
\end{Lemma}

We recall the following well-known $q$-product identity (for proof, see page 6 of \cite{Fine}).

\begin{Lemma}\label{lemma2}
Using the above notations, we have that
$$
\frac{(q;q)^2_{\infty}}{(q^2;q^2)_{\infty}}= 1 + 2 \sum_{n=1}^{\infty} (-1)^{n} q^{n^2}.
$$
\end{Lemma}

The following elementary lemma plays a crucial role in this paper.

\begin{Lemma} \label{lemma3}
If $\S$ is a subset of the positive integers, then the following are true:
\begin{displaymath}
F_{\S}(q)= \sum_{n \in \S} q^{n} \prod_{m=1}^{\infty} (1-q^{m+n})\\
= (q;q)_{\infty}\cdot  \sum_{\substack{ \lambda \in \mathcal{P} \\ \rm{lg}(\lambda) \in \S}} q^{\vert \lambda \vert}.
\end{displaymath}
\end{Lemma}
\begin{remark}
Lemma \ref{lemma3} may be viewed as a partition-theoretic case of Alladi's duality principle, originally stated in \cite{A} as a relation between functions on smallest versus largest prime divisors of integers, which he gave in full generality in a lecture \cite{Alladi_lecture}. 
\end{remark}

\begin{proof}
By inspection, we see that
$$
F_{\mathcal{S}}(q)=\sum_{\substack{ \lambda\in \mathcal P\\
{\text {\rm sm}}(\lambda)\in \mathcal{S}}} \mu^*_{\mathcal{P}}(\lambda) q^{|\lambda|}=
\sum_{n\in \mathcal{S}}q^n\prod_{m=1}^{\infty}(1-q^{m+n}).
$$
By factoring out $(q;q)_{\infty}$ from each summand, we find that
\begin{align*}
F_{\S}(q) &= \sum_{n \in \S} q^{n} \prod_{m=1}^{\infty} (1-q^{m+n}) 
= (q;q)_{\infty} \cdot \sum_{n \in \S} \frac{q^n}{(q;q)_{n}} \\
&= (q;q)_{\infty} \cdot \sum_{\substack{ \lambda \in \mathcal{P} \\ \rm{lg}(\lambda) \in \S}} q^{\vert \lambda \vert}.
\end{align*}
\end{proof}

\subsection{Case of $F_{S_{r,t}}(q)$}

Here we specialize Lemma~\ref{lemma3} to the sets  $\mathcal{S}_{r,t}$. The next lemma describes the $q$-series
$F_{\mathcal{S}_{r,t}}(q)$ in terms of a finite sum of quotients of infinite products. To prove this lemma we make
use of the $q$-Binomial Theorem.

\begin{Lemma} \label{finite}
If $t$ is a positive integer and $\zeta_t:=e^{2\pi i/t}$, then
\[F_{\S_{r,t}}(q) = (q;q)_{\infty} \cdot \frac{1}{t} \left[ \sum_{m=1}^{t} \frac{\zeta_{t}^{-mr}}{(\zeta_{t}^{m} q;q)_{\infty}} \right].\]
\end{Lemma}

\begin{proof}
From Lemma \ref{lemma3}, we have that
\[F_{\S_{r,t}}(q) = (q;q)_{\infty}\cdot \sum_{n=0}^{\infty} \frac{q^{tn+r}}{(q;q)_{tn+r}}.\]
By applying the $q$-Binomial Theorem (see Lemma \ref{lemma1}) with $a=0$ and $z=\zeta_t^m q$,  we find that
\begin{align*}
\frac{1}{t} \left[ \sum_{m=1}^{t} \frac{\zeta_{t}^{-mr}}{(\zeta_{t}^{m} q;q)_{\infty}} \right] = \frac{1}{t} \left[ \sum_{m=1}^{t} \sum_{n =0}^{\infty} \frac{ \zeta_{t}^{m(n-r)}q^{n}}{(q;q)_{n}} \right]. \end{align*}
Due to the orthogonality of roots of unity we have 
$$
\sum_{m=1}^{t} \zeta_{t}^{m(n-r)} = \begin{cases} t & \text{if} \ n \equiv r \pmod{t} \\
0 & \text{otherwise}. \end{cases}
$$
Hence, this sum allows us to sieve on the sum in $n$ leaving only those summands with $n\equiv r\pmod t$, namely
the series
\[\sum_{n =0}^{\infty} \frac{q^{tn+r}}{(q;q)_{tn+r}}.\]
Therefore, it follows that
\[F_{\S_{r,t}}(q) = (q;q)_{\infty} \cdot \frac{1}{t} \left[ \sum_{m=1}^{t} \sum_{n =0}^{\infty} \frac{ \zeta_{t}^{m(n-r)}q^{n}}{(q;q)_{n}} \right].\]
\end{proof}

\begin{Lemma} \label{limit}
If $a$ and $m$ are positive integers and $\zeta$   is a primitive $m$th root of unity, then
\[\lim_{q \to 1} \frac{(q;q)_{\infty}}{(\zeta^{a}q;q)_{\infty}} = \begin{cases} 1 & {\rm{if}} \ m \mid a \\
0 & \rm{otherwise}. \end{cases}\] 
\end{Lemma}

\begin{proof}
Since $(aq;q)_{\infty}^{\pm 1}$ is an analytic function of $q$ inside the unit disk (i.e., of $q:=e^{2\pi i z}$ with $z$ in the upper half-plane) when $|a|\leq 1$, the quotient on the left-hand side of Lemma \ref{limit} is well-defined as a function of $q$ (resp. of $z$), and we can take limits from inside the unit disk.  When $m \mid a$, the $q$-Pochhammer symbols cancel and the quotient is identically $1$.  When $m \nmid a$, then $(q;q)_{\infty}$ clearly vanishes as $q\to 1$ while $(\zeta^{a}q;q)_{\infty}$ is non-zero; thus the quotient is zero.
\end{proof} 

\section{Proofs of our results}

\subsection{Proof of Theorem \ref{12}}
Here we prove  the first part of Theorem \ref{12}; we defer the proof of the second
part until the next section because it is a special case of Theorem \ref{rt}.
\begin{proof} Here we prove the first part of Theorem~\ref{12}.
By Lemma \ref{finite} we have 
\begin{align*}
F_{\S_{1,2}}(q) &= (q;q)_{\infty} \cdot \frac{1}{2} \left[ \frac{1}{(q;q)_{\infty}} - \frac{1}{(-q;q)_{\infty}} \right] \\
&= \frac{1}{2} \left[ 1 - \frac{(q;q)_{\infty}}{(-q;q)_{\infty}} \right]\\
&=\frac{1}{2}\left [1-\frac{(q;q)_{\infty}^2}{(q^2;q^2)_{\infty}}\right ].
\end{align*}
Lemma \ref{lemma2} now implies that
\[F_{\S_{1,2}}(q) = \sum_{n=1}^{\infty} (-1)^{n+1}q^{n^2}.\]
To prove the $F_{\S_{2,2}}(q)$ identity, first recall that $\sum_{\lambda \in \mathcal{P}}\mu_{\mathcal{P}}^{*}(\lambda) q^{\vert \lambda \vert} = -(q;q)_{\infty}$.  From this we know $F_{\S_{1,2}}(q) + F_{\S_{2,2}}(q) = 1 - (q;q)_{\infty}$.  Using the identity for $F_{\S_{1,2}}(q)$ and Euler's Pentagonal Number Theorem completes the proof.
\end{proof}

\begin{proof} Now we prove Corollary 1.2.

Part 1. This claim follows easily from the first part of Theorem 1.1. The reader should recall
that $F_{S_{1,2}}(q)$ is the generating function for $\mu_{\mathcal{P}}^*(\lambda)=
-\mu_{\mathcal{P}}(\lambda)$.

Part 2. This corollary is not as immediate as the first part. Of course, we must classify the integer pairs $m$ and $n$ for which $n^2=m(3m-1)/2$. After simple manipulation, we find that this holds if and only if
$$
(6m-1)^2-6(2n)^2=1.
$$
In other words, we require that $(x,y)=(6m-1,2n)$ be a solution to Pell's equation
$$
x^2-6y^2=1.
$$
It is well known that all of the positive solutions to Pell's equation are of the form
$(x_k, y_k)$, where
$$
x_k +\sqrt{6}\cdot y_k = (5+2\sqrt{6})^k.
$$
Using this description and the elementary congruence properties of $(x_k,y_k)$, one easily obtains the second part of
Corollary 1.2.

\end{proof}

\subsection{Proof of Theorem \ref{rt}}
Here we prove the general limit formulas for the arithmetic densities of $\mathcal{S}_{r,t}$.
\begin{proof} Now we prove Theorem~\ref{rt}.
From Lemma \ref{finite} we have 
\[F_{\S_{r,t}}(q) = (q;q)_{\infty} \cdot \frac{1}{t} \left[ \sum_{m=1}^{t} \frac{\zeta_{t}^{-mr}}{(\zeta_{t}^{m} q;q)_{\infty}} \right].\]
We stress that we can take a limit here because we have a finite sum of functions which are analytic inside the unit disk.  Using Lemma \ref{limit} we see that 
\[\lim_{q \to 1} \frac{(q;q)_{\infty}}{(\zeta_{t}^{m} q;q)_{\infty}} = \begin{cases} 1 & {\rm{if}} \  m=t \\
0 & \rm{otherwise}. \end{cases}\]
From this we have 
\[\lim_{q\to 1} F_{\S_{r,t}}(q) = \frac{1}{t}.\]
The proof for $q \to \zeta$ where $\zeta$ is a primitive $m$th root of unity with $\gcd(m, t)=1$ follows the exact same steps.

\end{proof}

\subsection{Proofs of Theorem \ref{free} and Corollary \ref{free2}}
Here we will prove Theorem \ref{free} and Corollary \ref{free2} by building up $k$th power-free sets using arithmetic progressions.  We prove Theorem \ref{free} first.
\begin{proof} Here we prove Theorem \ref{free}.
The set of numbers not divisible by $p^{k}$ for any prime $p \leq N$ can be built as a union of sets of arithmetic progressions.  Therefore, for a given fixed $N$ we only need to understand divisibility by $p^{k}$ for all primes $p \leq N$.  Because the divisibility condition for each prime is independent from the other primes, we have
$$
F_{\S_{\rm{fr}}^{(k)}(N)}(q) = \sum_{\substack{0\leq r < M\\ p^k\nmid \ r}}
F_{\mathcal{S}_{r,M}}(q),
$$
where $M:=\prod_{\substack{p\leq N \ \\ {\text {\rm prime}}}} p^k$.
We have a finite number of summands, and the result now follows immediately from
Theorem~\ref{rt}.
\end{proof}

Next, we will prove Corollary \ref{free2}.
\begin{proof} Here we prove Corollary \ref{free2}.
For fixed $N$ define 
$$\zeta_{N}(k) := \prod_{p \leq N \rm{prime}} \left(\frac{1}{1-p^k} \right),
$$ so $\lim_{q \to 1} F_{\S_{\rm{fr}}^{(k)}(N)}(q)$ $= \frac{1}{\zeta_{N}(k)}$.  It is clear $\lim_{N \to \infty} \zeta_{N}(k)= \zeta(k)$.  It is in this sense that we say $\lim_{q \to 1} F_{\S_{\rm{fr}}^{(k)}}(q) = \frac{1}{\zeta(k)}$.
\end{proof}

\section{Examples}

\begin{Example}
In the case of $\S_{1,3}$, which has arithmetic density 1/3, Theorem \ref{rt} holds for any $m$th root of unity where $3 \nmid m$.  The two tables below illustrate this as $q$ approaches $\zeta_{1} = 1$ and $\zeta_{4} = i$, respectively.

\begin{center}
\begin{tabular}{ | c | c | }
\hline
$q$ & $F_{\mathcal S_{1,3}}(q)$\\ \hline
$0.70$ & $0.340411885...$ \\ \hline
$0.75$ & $0.335336994...$ \\ \hline
$0.80$ & $0.333552814...$ \\ \hline
$0.85$ & $0.333331545...$ \\ \hline
$0.90$ & $0.333333329...$ \\ \hline
$0.95$ & $0.333333333...$ \\   
\hline
\end{tabular}
\end{center}
\vspace{0.5cm}


\begin{center}
\begin{tabular}{ | c | c | }
\hline
$q$ & $F_{\mathcal S_{1,3}}(q)  $\\ \hline
$0.70i$ & $\approx 0.034621+0.793781i$ \\ \hline
$0.75i$ & $\approx 0.057890+0.802405i$ \\ \hline
$0.80i$ & $\approx 0.097030+0.771774i$ \\ \hline
$0.85i$ & $\approx 0.167321+0.674712i$ \\ \hline
$0.90i$ & $\approx 0.294214+0.454400i$ \\ \hline
$0.95i$ & $\approx 0.424978+0.067775i$ \\ \hline 
$0.97i$ & $\approx 0.376778-0.016187i$ \\ \hline
$0.98i$ & $\approx 0.340170+0.005772i$ \\ \hline
$0.99i$ & $\approx 0.332849+0.000477i$ \\  
\hline
\end{tabular}
\end{center}
\vspace{0.5cm}
\end{Example}

\begin{Example}
The table below illustrates Theorem \ref{free} for the set $\mathcal S^{(2)}_{\rm{fr}}(5)$, which has arithmetic density $16/25 = 0.64$.  These are the positive integers which are not divisible by 4, 9 and 25.
Here we give numerics for the case of $F_{\mathcal S^{(2)}_{\rm{fr}}(5)} (q)$ as $q\to 1$ along the real axis.
 
\begin{center}
\begin{tabular}{ | c | c | }
\hline
$q$ & $F_{\mathcal S^{(2)}_{\rm{fr}}(5)} (q)$\\ \hline
$0.90$ & $0.615367...$ \\ \hline
$0.91$ & $0.619346...$ \\ \hline
$0.92$ & $0.625991...$ \\ \hline
$0.93$ & $0.631607...$ \\ \hline
$0.94$ & $0.631748...$ \\ \hline
$0.95$ & $0.631029...$ \\ \hline
$0.96$ & $0.638291...$ \\ \hline
$0.97$ & $0.639893...$ \\  
\hline
\end{tabular}
\end{center}
\vspace{0.5cm}

\end{Example}


\begin{Example}
Here we approximate the density of $\mathcal S_{\rm{fr}}^{(4)}$, the fourth power-free positive integers. Since $\zeta(4)=\pi^4/90$, it follows that the arithmetic density of $\mathcal S_{\rm{fr}}^{(4)}$ is
${90}/{\pi^4} \approx 0.923938...$. Here we choose $N=5$ and compute the arithmetic density of $\mathcal S_{\rm{fr}}^{(4)}(5)$, the positive integers which are not divisible by $2^4, 3^4$, and $5^4$. The density
of this set is $208/225 \approx 0.924444...$. This density is fairly close to the density of fourth power-free integers because the convergence in the $N$ aspect is significantly faster for fourth power-free integers than for square-free integers, as discussed above.
 
\begin{center}
\begin{tabular}{ | c | c | }
\hline
$q$ & $F_{\mathcal S_{\rm{fr}}^{(4)}(5)}(q)$\\ \hline
$0.90$ & $0.934926...$ \\ \hline
$0.91$ & $0.936419...$ \\ \hline
$0.92$ & $0.936718...$ \\ \hline
$0.93$ & $0.935027...$ \\ \hline
$0.94$ & $0.931517...$ \\ \hline
$0.95$ & $0.925619...$ \\ \hline
$0.96$ & $0.921062...$ \\ \hline
$0.97$ & $0.925998...$ \\ \hline
$0.98$ & $0.924967...$ \\  
\hline
\end{tabular}
\end{center}
\vspace{0.5cm}
\end{Example}

\end{document}